\documentclass{article}
\usepackage{amsmath, amsthm, amssymb, fullpage, url}

\newtheorem{theorem}{Theorem}
\newtheorem{lemma}[theorem]{Lemma}

\title{Farey Statistics in Time $n^{2/3}$ and \\
  Counting Primitive Lattice Points in Polygons}
\author{Mihai P\v{a}tra\c{s}cu \\ {\tt mip@mit.edu}}

\newcommand{\Ot}{\widetilde{O}}
\newcommand{\Far}{\mathcal{F}}

\newcommand{\Stat}{\mbox{\sc Statistic}}
\newcommand{\Rank}{\mbox{\sc Rank}}

\newcommand{\twodots}{\mathinner{\ldotp\ldotp}}
\newcommand{\poly}{\mathop{\rm poly}\nolimits}

\begin{document}

\maketitle

\begin{abstract}
We present algorithms for computing ranks and order statistics in the
Farey sequence, taking time $\Ot(n^{2/3})$. This improves on the
recent algorithms of Pawlewicz~\cite{pawlewicz07farey}, running in
time $\Ot(n^{3/4})$. We also initiate the study of a more general
algorithmic problem: counting primitive lattice points in planar
shapes.
\end{abstract}

\noindent
{\bf Since the publication of this technical report, this work has
been extended and merged with the paper of Pawlewicz. The merged
version is available at:}\\
\url{http://web.mit.edu/~mip/www/papers/farey2/paper.pdf}

\section{An Improved Algorithm for the Farey Sequence}

The Farey sequence of order $n$, denoted $\Far_n$, is the ordered list
of irreducible fractions $\frac{a}{b}$ with $a \le b \le n$. This
sequence is a well-studied mathematical object, with fascinating
properties. See \cite{gkp-concretebook} for a discussion at length.
The sequence has $\Theta(n^2)$ terms, and there are many algorithms
for generating it entirely in $O(n^2)$ time. Perhaps the best known
ones are based on the Stern-Brocot tree, and the properties of the
mediant.

One can ask, however, for more local access to the sequence. Two
natural questions arise:

\begin{itemize}
\item given a number $x \in [0,1]$ find $\Rank(x,n) = \big| \Far_n
  \cap [0,x] \big|$.

\item given an index $k \le |\Far_n|$ find $\Stat(k,n) =$ the $k$-th
  value in $\Far_n$ (in sorted order). 
\end{itemize}

The $\Stat$ problem can be solved with $O(\lg n)$ calls to the $\Rank$
problem~\cite{pawlewicz07farey}, so below only bounds for the $\Rank$
version are discussed.

To the best of my knowledge, the question was first formulated in
2003, when I proposed it as a contest problem at the 11th Balkan
Olympiad in Informatics. The official solution consisted of an
$O(n\lg n)$ algorithm, and several contestant also found this
algorithm. We~\cite{patrascu04farey} later described a solution with a
slightly better running time of $O(n)$.

Quite recently, Pawlewicz~\cite{pawlewicz07farey} broke the linear
time barrier, and provided an algorithm with running time
$O(n^{3/4})$. In the present, I describe an improved $O(n^{2/3}
\lg^{1/3} n)$ algorithm.

\subsection{Review: The Algorithm of Pawlewicz}

Let $S_n(x) = \left| \left\{ \frac{a}{b} \mid b \le n ~\land~
\frac{a}{b} \le x ~\land~ \gcd(a,b) = 1 \right\} \right|$. This
provides the quantity we want to compute. It can be seen that:
\begin{equation}
S_n(x) = \sum_{b=1}^n \lfloor bx \rfloor - \sum_{d\ge 2}
  S_{\lfloor \frac{n}{d} \rfloor}(x)       \label{eq:my_rec}
\end{equation}
It is shown in \cite{pawlewicz07farey} that $A_n(x) = \sum_{b=1}^n
\lfloor bx \rfloor$ can be computed in $O(\lg n)$ time. 

Thus, the only challenge is to estimate the recursive component of the
sum. The recursion will only need $S_{\lfloor n/d \rfloor}(x)$ for all
$d$, since $\big\lfloor \lfloor \frac{n}{d_1} \rfloor / d_2
\big\rfloor = \big\lfloor \frac{n}{d_1 d_2} \big\rfloor$. 

The crux of the algorithm is the following observation: given all
relevant $S_i(x)$, for $i < k$, then $S_k(x)$ can be computed in
$O(\sqrt{k})$ time. Indeed, $\sum_{d\ge 2} S_{\lfloor k/d \rfloor}(x)$
only contains at most $2 \sqrt{k}$ distinct terms: $\sqrt{k}$ terms
corresponding to $d \le \sqrt{k}$, and at most $\sqrt{k}$ terms for $d
\ge \sqrt{k}$ because then we have $k/d \le \sqrt{k}$. The latter
terms have multiplicities, but the multiplicity of each term can be
computed easily in $O(1)$ time.

Applying this observation to compute all needed $S$ terms
recursively, the running time is:
\[ \sum_{d=1}^n \sqrt{\frac{n}{d}} ~\le~
  \sum_{d=1}^{\sqrt{n}} \sqrt{\frac{n}{d}} + \sum_{j=1}^{\sqrt{n}} \sqrt{j}
~\le~ \sqrt{n} \cdot \sum_{d=1}^{\sqrt{n}} \frac{1}{\sqrt{d}}
  + \sum_{j=1}^{\sqrt{n}} \sqrt{\sqrt{n}}
~\le~ \sqrt{n}\cdot O(\sqrt[4]{n}) + \sqrt{n} \cdot \sqrt[4]{n}
~=~  O(n^{3/4})
\]

\subsection{Our Improved Algorithm}

The key to our improved algorithm is to show that $S_1, \dots, S_k$
can be computed in time $O(k\lg k)$. Then, the remaining terms can be
computed by the old algorithm in time $\sum_{d=1}^{n/k} \sqrt{n/d} =
\sqrt{n} \cdot O(\sqrt{n/k}) = O(n/\sqrt{k})$. The total running time
is then $O(k\lg k + n/\sqrt{k})$, so it is optimized by picking $k =
(n/\lg n)^{2/3}$. Thus, the running time is $O(n^{2/3} \lg^{1/3} n)$.

To compute $S_1, \dots, S_k$ efficiently, we make the observation that
the composition (with respect to recursive terms) of $S_i$ and
$S_{i-1}$ are not too different. Specifically, $\sum_{d\ge 2}
S_{\lfloor i/d \rfloor}(x)$ differs from $\sum_{d\ge 2} S_{\lfloor
(i-1)/d \rfloor}(x)$ only for the values of $d$ that $i$ is a multiple
of.

To maintain understanding of divisibility by all $d$'s, and compute
$S_i(x)$ values in order, we use an algorithm similar in flavor to
Eratosthene's sieve \cite{eratosthene}. We first create an array
$D[1\twodots k]$, where $D[i]$ holds a list of all divisors of $i$. To
create the array, simple consider all $d$, and add $d$ to $D[md]$, for
all $m$. This takes time $\sum_{d\ge 1} \frac{k}{d} = O(k\lg k)$.

Now, iterate $i$ from $1$ to $k$, maintaining $\sum_{d\ge 2}
S_{\lfloor i/d \rfloor}(x)$ at all times. The sum is updated by
considering all $d\in D[i]$, subtracting $S_{\lfloor (i-1)/d \rfloor}(x)$
and adding $S_{\lfloor i/d \rfloor}(x)$. The complexity is linear in
the size of $D[1\twodots k]$, and is thus $O(k\lg k)$. To compute
$S_i(x)$ from this running sum, we only need $A_i(x)$, which takes
$O(\lg i)$, giving an additive $O(k\lg k)$.

\section{Counting Primitive Lattice Points}

A primitive lattice point is a point $(x,y)$ in the plane, with $x, y
\in \mathbb{Z}$ and $\gcd(x,y) = 1$. We observe that computing
$\Rank(x,n)$ in the Farey sequence is equivalent to couting primitive
lattice points inside the right triangle defined by $(0,0)$, $(n,0)$
and $(n, xn)$. Let us now generalize the algorithm to counting
primitive lattice points in more general shapes.

Counting primitive lattice points in planar shape is a relatively new
topic in mathematics, but one that is gathering significant momentum
\cite{moroz85primitive, nowak88primitive, hensley94primitive,
nowak95primitive, huxley96primitive, muller96primitive,
nowak97primitive, kratzel99primitive, zhai99primitive,
boca00primitive, wu02primitive, zhai03primitive, nowak05primitive}.
In the mathematical sense, ``counting'' refers to estimating the
number of primitive points with a small error, as the size of the
shape goes to infinity.

In this paper, we initiate the study of the \emph{algorithmic} problem
of counting (exactly) the number of primitive lattice points inside a
given shape. More precisely, we study this problem for polygons
containing the origin. The condition that the shape should contain the
origin also appears in the mathematical works referenced above, and is
natural given that we are counting points visible from the origin.

\begin{theorem}  \label{thm:primitive}
Let $P$ be a polygon containing the origin, defined by $k$ vertices at
$b$-bit rational coordinates. If $D$ is the diameter of the polygon,
one can count the number of primitive lattice points inside $P$ in
time $D^{6/7} \cdot k \cdot b^{O(1)}$.
\end{theorem}

A very pertinent question is how efficient this running time actually
is. Remember that \cite{patrascu04farey} shows that the Farey rank problem
can be used to factor integers. Since counting primitive lattice
points is a generalization, we conclude that a polynomial-time
algorithm, i.e.~$\poly(k\cdot b)$, is likely impossible. Thus, the
algorithm needs to depend on some parameter describing the polygon,
which can be exponential in $b$. One such parameter is the diameter
$D$. Clearly, however, this is not the only choice, and it is
conceivable that other measures lead to better results. For right
triangles such as those in the Farey rank problem, the diameter is
$n$, yet we know a better algorithm with time essentially $n^{2/3}$.

A trivial alternative to Theorem~\ref{thm:primitive} is the algorithm
which iterates over all lattice points inside the polygon, and runs
Euclid's algorithm on each point. It is possible 
\cite{barvinok94lattice,beck02lattice,hiroki05lattice}
to list all lattice points with a polynomial $\poly(k\cdot b)$ cost
per point. If the polygon has integral coordinates, Pick's formula
shows that the number of lattice points inside is asymptotically equal
to the area. Thus, the exhaustive algorithm has complexity
proportional to the area, times polynomial factors.

Unfortunately, the area and the diameter are not related in the
worst-case (e.g., for very skinny shapes).  However, in the more
``typical'' case when the polygon is fat, the area is
$A=\Theta(D^2)$. Thus our running time of $D^{6/7}$ can be rewritten
as $A^{3/7}$, which gives a significant saving over the exhaustive
algorithm. It is an interesting open problem to construct an algorithm
which beats exhaustive search for \emph{any} polygon.

\paragraph{A standard idea.}
Let $P$ be a polygon, defined by rational points $(x_1, y_1), \dots,
(x_k, y_k)$. Then, let $P_{/d}$ be the polygon defined by $(\frac{x_1}{d},
\frac{y_1}{d}), \dots, (\frac{x_n}{d}, \frac{y_n}{d})$. Define $A(P)$
to be the number of lattice points inside polygon $P$, and $S(P)$ the
number of primitive lattice points inside $P$.

We first observe the following recursive formula:
\begin{equation}
  S(P) = A(P) - \sum_{d \ge 2} S(P_{/d})   \label{eqn:polyrec}
\end{equation}
Indeed, every point in $A(P)$, but not in $S(P)$ is a nonprimitive
lattice point $(x,y)$. If $\gcd(x,y) = d > 1$, then $(\frac{x}{d},
\frac{y}{d})$ is a primitive lattice point. Furthermore, such a point
is inside $P_{/d}$, so we can remove all points with a greater common
divisor equal to $d$ by subtracting $P_{/d}$.

We can bound the recursion depth in \eqref{eqn:polyrec}, by appealing
to the diameter $D$. We note that the diameter of $P_{/(D+1)}$ is less
than 1, so it does not contain any lattice point outside the
origin. Thus, $S(P_{/d}) = 0$ for all $D < d < \infty$, and
$S(P_{/\infty}) = 1$ (the origin). Then, it suffices to consider only
$P, P_{/2}, \dots, P_{/D}$ in the algorithm.

To compute the ``constants'' $A(P_{/i})$ in the recursion, one needs
to compute the number of lattice points inside polygons with rational
coordinates. This is a well-studied problem, and 
\cite{barvinok94lattice,beck02lattice,hiroki05lattice}
given polynomial-time algorithms.  Observe that for $i\le D \le 2^b$,
coordinates of $P_{/i}$ have at most $2b$ bits of precision.  Thus,
computing $A(P_{/i})$ takes time $O(k \cdot \poly(b))$; the dependence
on $k$ is linear by triangulating the polygon.

Note that formula \eqref{eqn:polyrec} is very similar to the recursive
formula for the Farey rank problem \eqref{eq:my_rec}. Indeed,
\eqref{eq:my_rec} is simply a transcription of \eqref{eqn:polyrec},
where the polygons are the relevant right triangles.

It is thus tempting to conjecture that we can use the same dynamic
program to evaluate \eqref{eqn:polyrec}. Unfortunately, this is not
the case, for a somewhat subtle reason. Due to the geometry of the
rank problem, $S_{n/d}$ was the same as $S_{\lfloor n/d \rfloor}$. By
taking floors, we concluded that among $S_{\lfloor n/d \rfloor}$'s with
$d\in \{ \sqrt{n}, \dots, n\}$, there are only $\sqrt{n}$ distinct
quantities. Unfortunately, in general we cannot round the vertices of
$P_{/d}$ to lattice points, and thus we cannot conclude that for $d
\ge \sqrt{D}$ there are only $\sqrt{D}$ distinct cases.

\newcommand{\Os}{O^*}

\paragraph{A more careful analysis.}
Since we are not interested in factors of $k \cdot b^{O(1)}$ in the
running time, let us define $\Os(f) = f \cdot k \cdot b^{O(1)}$.
To use ideas similar to the previous dynamic program, we begin by
obtaining a weaker bound for computing the small terms of the
recurrence:

\begin{lemma}  \label{lem:precomp}
We can (implicitly) compute $S(P_{/\tau}), S(P_{/(\tau+1)}), \dots,
S(P_{/D})$ in time $\Os\big( \frac{D^2}{\tau^2} \big)$.
\end{lemma}

\begin{proof}
Note that $P_{/D} \subseteq P_{/(D-1)} \subseteq \cdots \subseteq
P_{/\tau}$. Also, the diameter of $P_{/\tau}$ is $D/\tau$, implying that
$P_{/\tau}$, and any smaller polygon only contain $O\big( (\frac{D}{\tau})^2
\big)$ lattice points. Since $S(P_{/i})$ can only grow between $1$ and
$O\big( (\frac{D}{\tau})^2 \big)$ when $i$ goes from $D$ to $\tau$, there
are only so many distinct values that can appear.

To actually compute these values, we perform an exhaustive enumeration
of all lattice points in $P_{/\tau}$. Points which are not primitive are
discarded. For every primitive point $(x,y)$, we compute a value
$\varphi(x,y)$ which is the minimum $i$ such that $P_{/i}$ does not
contain it. This is done by a binary search for $i$. Every comparison
is a point-in-polygon test, which takes $O(k)$ time.

We now sort the $\phi(x,y)$ values from largest to smallest. The
sorted list gives us an \emph{implicit} representation for $S(P_{/i})$
for every $i \ge \tau$. Indeed, it suffices to binary search for the
first occurrence of $i$; the elements to the left correspond to
primitive points which are inside $P_{/i}$.
\end{proof}

Let us now consider the problem of computing a term using
our recursion: $S(P_{/i}) = A(P_{/i}) - \sum_{d\ge 2} S(P_{/id})$.
We assume terms $S(P_{/j})$ for $j \ge i$ have already been computed;
in particular, for $j \ge \tau$ we have the implicit representation
from Lemma~\ref{lem:precomp}.

In principle, the expression for $S(P_{/i})$ has $\lfloor \frac{D}{i}
\rfloor$ terms. However, as in Lemma~\ref{lem:precomp}, we can observe
that for fixed $\delta$, there are only $O\big( \frac{D^2}{\delta^2}
\big)$ distinct terms among all $S(P_{/id})$'s with $id \ge
\delta$. These terms can actually be summed up in $\Os \big(
\frac{D^2}{\delta^2} \big)$ time. Indeed, we begin with $d_0 = \lceil
\frac{\delta}{i} \rceil$, and binary search for the minimum $d_1$ such
that $S(P_{/id_1}) < S(P_{/id_0})$. We add to the sum $S(P_{/id_0})
\cdot (d_1 - d_0)$, then binary search for the minimum $d_2$ leading
to a different value, and so on.

After dealing like this with all terms $d \ge \frac{\delta}{i}$, we
can simply add the remaining terms. The running time for computing
$S(P_{/i})$ is therefore $\Os \big( \frac{D^2}{\delta^2} +
\frac{\delta}{i} \big)$. We can optimize by setting $\delta = D^{2/3}
i^{1/3}$, yielding a running time of $\Os \big( (\frac{D}{i} )^{2/3}
\big)$.

We now want to optimize the parameter $\tau$ in Lemma~\ref{lem:precomp}.
A choice of $\tau$ implies that we spend $\Os\big( \frac{D^2}{\tau^2}
\big)$ time by Lemma~\ref{lem:precomp}, and then run the dynamic
program for terms $S(P_{/i})$ with $i < \tau$. This second part will take time:
\[
  \sum_{i=1}^{\tau} \Os \bigg( \frac{D}{i} \bigg)^{2/3}
= \Os( D^{2/3} \tau^{1/3})
\]
Then, we can optimize by setting $\frac{D^2}{\tau^2} = D^{2/3}
\tau^{1/3}$, i.e.~$\tau = D^{4/7}$. The final running time is
therefore $\Os( D^{6/7} )$.

\bibliographystyle{alpha}
\bibliography{../../general,../primitive}

\end{document}